\documentstyle [12pt,epsf, graphicx, amssymb]{article}

\begin{document}
\def\l{\lambda}
\def\m{\mu}
\def\a{\alpha}
\def\b{\beta}
\def\g{\gamma}
\def\d{\delta}
\def\e{\epsilon}
\def\o{\omega}
\def\O{\Omega}
\def\v{\varphi}
\def\t{\theta}
\def\r{\rho}
\def\bs{$\blacksquare$}
\def\bp{\begin{proposition}}
\def\ep{\end{proposition}}
\def\bt{\begin{th}}
\def\et{\end{th}}
\def\be{\begin{equation}}
\def\ee{\end{equation}}
\def\bl{\begin{lemma}}
\def\el{\end{lemma}}
\def\bc{\begin{corollary}}
\def\ec{\end{corollary}}
\def\pr{\noindent{\bf Proof: }}
\def\note{\noindent{\bf Note. }}
\def\bd{\begin{definition}}
\def\ed{\end{definition}}
\def\C{{\mathbb C}}
\def\P{{\mathbb P}}
\def\Z{{\mathbb Z}}
\def\d{{\rm d}}
\def\deg{{\rm deg\,}}
\def\deg{{\rm deg\,}}
\def\arg{{\rm arg\,}}
\def\min{{\rm min\,}}
\def\max{{\rm max\,}}

\newcommand{\norm}[1]{\left\Vert#1\right\Vert}
\newcommand{\abs}[1]{\left\vert#1\right\vert}

\newcommand{\set}[1]{\left\{#1\right\}}
\newcommand{\setb}[2]{ \left\{#1 \ \Big| \ #2 \right\} }

\newcommand{\IP}[1]{\left<#1\right>}
\newcommand{\Bracket}[1]{\left[#1\right]}
\newcommand{\Soger}[1]{\left(#1\right)}

\newcommand{\Integer}{\mathbb{Z}}
\newcommand{\Rational}{\mathbb{Q}}
\newcommand{\Real}{\mathbb{R}}
\newcommand{\Complex}{\mathbb{C}}

\newcommand{\eps}{\varepsilon}
\newcommand{\To}{\longrightarrow}
\newcommand{\varchi}{\raisebox{2pt}{$\chi$}}

\newcommand{\E}{\mathbf{E}}
\newcommand{\Var}{\mathrm{var}}

\def\squareforqed{\hbox{\rlap{$\sqcap$}$\sqcup$}}
\def\qed{\ifmmode\squareforqed\else{\unskip\nobreak\hfil
\penalty50\hskip1em\null\nobreak\hfil\squareforqed
\parfillskip=0pt\finalhyphendemerits=0\endgraf}\fi}

\renewcommand{\th}{^{\mathrm{th}}}
\newcommand{\Dif}{\mathrm{D_{if}}}
\newcommand{\Difp}{\mathrm{D^p_{if}}}
\newcommand{\GHF}{\mathrm{G_{HF}}}
\newcommand{\GHFP}{\mathrm{G^p_{HF}}}
\newcommand{\f}{\mathrm{f}}
\newcommand{\fgh}{\mathrm{f_{gh}}}
\newcommand{\T}{\mathrm{T}}
\newcommand{\K}{^\mathrm{K}}
\newcommand{\PghK}{\mathrm{P^K_{f_{gh}}}}
\newcommand{\Dig}{\mathrm{D_{ig}}}
\newcommand{\for}{\mathrm{for}}
\newcommand{\End}{\mathrm{end}}

\newtheorem{th}{Theorem}[section]
\newtheorem{lemma}{Lemma}[section]
\newtheorem{definition}{Definition}[section]
\newtheorem{corollary}{Corollary}[section]
\newtheorem{proposition}{Proposition}[section]

\begin{titlepage}

\begin{center}

\topskip 5mm

{\LARGE{\bf {Higher derivatives of functions

\vskip 4mm

vanishing on a given set}}}

\vskip 8mm

{\large {\bf Y. Yomdin}}

\vspace{6 mm}

{Department of Mathematics, The Weizmann Institute of Science,
Rehovot 76100, Israel. e-mail: yosef.yomdin@weizmann.ac.il}

\end{center}

\vspace{6 mm}
\begin{center}

{ \bf Abstract}
\end{center}

{\small Let $f: B^n \rightarrow {\mathbb R}$ be a $d+1$ times continuously differentiable function on the unit ball $B^n$, with $\max_{z\in B^n} \Vert f(z) \Vert=1$. A well-known fact is that if $f$ vanishes on a set $Z\subset B^n$ with a non-empty interior, then for each $k=1,\ldots,d+1$ the norm of the $k$-th derivative $||f^{(k)}||$ is at least $M=M(n,k)>0$.

\medskip

We show that this fact remains valid for all ``sufficiently dense'' sets $Z$ (including finite ones). The density of $Z$ is measured via the behavior of the covering numbers of $Z$. In particular, the bound $||f^{(k)}||\ge \tilde M=\tilde M(n,k)>0$ holds for each $Z$ with the box (or Minkowski, or entropy) dimension $\dim_e(Z)$ greater than $n-\frac{1}{k}$.}

\end{titlepage}

\newpage


\section{Introduction}
\setcounter{equation}{0}


\smallskip


\smallskip

In this paper we study a very special setting of the Whitney smooth extension problem (\cite{Bru.Shv,Fef,Fef.Kla,Whi1,Whi2,Whi3}). Let $Z\subset B^n \subset {\mathbb R}^n$ be a closed subset of the unit ball $B^n$. We look for $C^{d+1}$-smooth functions $f:B^n\to {\mathbb R}$, vanishing on $Z$. Such $C^{d+1}$-smooth (and even $C^\infty$) functions $f$ always exist, since any closed set $Z$ is a set of zeroes of a $C^\infty$-smooth function.

\medskip

We normalize the extensions $f$ requiring $\max_{B^n}|f|=1$, and ask for the minimal possible norm of the last derivative $||f^{(d+1)}||$, which we call {\it the $d$-rigidity ${\cal RG}_d(Z)$ of $Z$}. In other words, for each normalized $C^{d+1}$-smooth function $f:B^n\to {\mathbb R},$ vanishing on $Z$, we have
$$
||f^{(d+1)}||\ge {\cal RG}_d(Z),
$$
and ${\cal RG}_d(Z)$ is the maximal number with this property.

\smallskip

Recent exciting developments in the general Whitney problem (see \cite{Bru.Shv,Fef,Fef.Kla} and references therein), provide essentially a complete answer to the general Whitney extension problem in any dimension. In particular, as in classical Whitney's results in dimension one (\cite{Whi2}), it is enough to check only finite subsets of $Z$ with cardinality bounded in terms of $n$ and $d$ only. There is also an algorithmic way to estimate the minimal extension norm for any finite $Z$. However, a possibility of an explicit answer, as in dimension one, through a kind of multi-dimensional divided finite differences, remains an open problem.

\smallskip

Of course, the results of \cite{Fef,Fef.Kla} provide, in particular, an algorithmic way to estimate ${\cal RG}_d(Z)$ for any finite $Z$. However, our goal in the present paper, as well as in our previous papers \cite{Yom1,Yom5,Yom6}, related to smooth rigidity, is somewhat different: we look for an explicit answer in terms of simple, and directly computable geometric (or topological) characteristics of $Z$.

\smallskip

Let us now state the main results of this paper. From now on we always assume that an integer $d\ge 1$ is fixed. In order to explain the question we deal with, first we shortly recall some of basic properties of the $d$-rigidity ${\cal RG}_d(Z)$. See Section \ref{Sec:background} below for a more detailed presentation.

\medskip

In dimension one, ${\cal RG}_d(Z)=0,$ if the cardinality $|Z|\le d$, and ${\cal RG}_d(Z)\ge \frac{(d+1)!}{2^{d+1}}$ otherwise. In higher dimensions ${\cal RG}_d(Z)$ may attain arbitrarily small positive values. However, {\it for $Z$ with a non-empty interior we always have ${\cal RG}_d(Z)\ge \frac{(d+1)!}{2^{d+1}}$, independently of the size and the geometry of $Z$.}

\medskip

In \cite{Yom5,Yom6} the following question was discussed: {\it Can this last property be extended to other $Z$, beyond those with a non-empty interior? In particular, is it true for sufficiently dense finite sets $Z$?}

\medskip

Here is a partial answer:

\bt\label{thm:main.intro1}
If the box dimension $dim_e(Z)$ is greater than $n-\frac{1}{d+1}$, then
$$
{\cal RG}_d(Z)\ge M=M(n,d)>0,
$$
where the positive constant $M$ depends only on $n$ and $d$.
\et
Definition of the box (or Minkowski, or entropy ...) dimension is given in Section \ref{Sec:proof.main.res} below. In particular, the result of Theorem \ref{thm:main.intro1} provides examples of discrete, but sufficiently dense, sets $Z$, for which ${\cal R}_d(Z)$ behaves in the same way as for sets with a non-empty interior.

\medskip

A non-asymptotic version of this result provides examples of finite, but sufficiently dense, sets $Z$ with the same property. It requires the following definition:

\bd\label{def:h.regular}

A finite set $Z$ inside the cube $Q^n_s=[0,s]^n$ of size $s$ is called $h$-dense, if the following condition is satisfied:

\smallskip

\noindent In each $h$-sub-cube of $Q^n_s$ there is a point of $Z$.

\ed
In particular, if we start with a regular $h/2$-grid in $Q^n_s$, and shift its points to distances at most $\frac{h}{8},$ we obtain an $h$-dense set.

\smallskip


\bt\label{main:rigid.h.regular}
Let $Z\subset Q^n_s$ be an $h$-dense set. Then if $h\le \xi s^{n(d+1)}$, then
$$
{\cal RG}_d(Z)\ge M>0,
$$
with the constant $\xi$ depending only on $n$ and $d$.
\et
Let us describe our basic approach. If we could find a straight line $\ell$ in ${\mathbb R}^n$, passing through the point $z_0$, where the absolute value $|f(z)|$ is equal to one, and through some $d+1$ distinct points in $Z$, we could immediately get the required lower bound ${\cal RG}_d(Z)\ge \frac{(d+1)!}{2^{d+1}}$ on the $d+1$-st derivative of $f$, via the basic properties of $d$-rigidity, mentioned above. However, for a generic finite set $Z$ any straight line $\ell$ meets $Z$ at one or two points, at most. Instead we replace $\ell$ by a smooth curve $\o$, and try to mimic the calculations for $\ell$. This requires analysis of the high order chain-rule expressions, on one side, and construction of curves $\o$ with small high-order derivatives, passing through some $d+1$ distinct points in $Z$, on the other side. We show, using a kind of ``discrete integral geometry'', that already for finite or discrete sets $Z$, which are dense enough (in particular, when $dim_e(Z) > n-\frac{1}{d+1}$), the required curves exist.

\medskip

The paper is organized as follows: in Section \ref{Sec:background} we recall some basic properties of $d$-rigidity, as well as some old and some very recent results of \cite{Yom1,Yom5,Yom6}. In Section \ref{Sec:Comparison.Deriv} we provide the first main technical ingredient of our approach - the analysis of the high order chain-rule expressions. In Section \ref{Sec:thickness} we define a certain geometric characteristic of sets $Z$, responsible for the existence of smooth curves $\o$ with ``small'' high-order derivatives, passing through at least $d+1$ points of $Z$. In Section \ref{sec:entr.bounds.nu} we provide the second main technical ingredient of our approach - a ``discrete integral geometry''. We show that for sufficiently dense sets $Z$ the required smooth curves $\o$ can be constructed. Finally, in Section \ref{Sec:proof.main.res}, we combine the tools of the previous sections in order to proof the main results of the paper.


\section{Some backgound}\label{Sec:background}
\setcounter{equation}{0}

Let $f: B^n \rightarrow {\mathbb R}$ be a $d+1$ times continuously differentiable function on $B^n$. For $l=0,1,\ldots,d+1$ put
$$
M_l(f)=\max_{z\in B^n} \Vert f^{(l)}(z) \Vert,
$$
where the pointwise norm $\Vert f^{(l)}(z) \Vert$ of the $l$-th derivative $f^{(l)}(z)$  of $f$ at $z$ is defined as the maximum of the absolute values of all the partial derivatives of $f$ of order $l$ at $z$.

\medskip

By technical reasons we consider below only closed sets $Z$, which are contained in a concentric ball $\hat B^n$ of radius $\frac{1}{3}$. For $Z\subset \hat B^n$ let $U_d(Z)$ denote the set of $C^{d+1}$ smooth functions $f(z)$ on $B^n$, vanishing on $Z$, with $M_0(f)=1$.

\bd\label{def:rigidity}
For $Z\subset \hat B^n$ we define the $d$-th rigidity constant ${\cal RG}_d(Z)$ as
$$
{\cal RG}_d(Z)=\inf_{f\in U_d(Z)}M_{d+1}(f).
$$
\ed
By this definition we get immediately $M_{d+1}(f)\ge {\cal RG}_d(Z)$ for any $f(z)$ on $B^n$, vanishing on $Z$, with $M_0(f)=1$. Our goal is to estimate ${\cal RG}_d(Z)$ in terms of accessible geometric features of $Z$.

\medskip

Notice that we do not insist on $Z$ being exactly the set of zeroes $Y(f)$ of the functions $f\in U_d(Z)$, but just require $Z\subset Y(f)$.

\medskip

As an example, consider the case of dimension $n=1$. Here we have the following important fact (see e.g. \cite{Yom5}):

\bp\label{prop:d.points}
For any $Z\subset \hat B^1$ we have ${\cal RG}_d(Z)\ge \frac{(d+1)!}{2^{d+1}},$ if $Z$ consists of at least $d+1$ different points, and ${\cal RG}_d(Z)=0$ if $Z$ consists of at most $d$ different points.
\ep

Thus in dimension one the minimal non-zero value of ${\cal RG}_d(Z)$ is $\frac{(d+1)!}{2^{d+1}}.$ This is not true any more in higher dimensions: for $Z\subset \hat B^n, \ n\ge 2,$ the $d$-rigidity ${\cal RG}_d(Z)$ attains arbitrarily small positive values (\cite{Yom6}. This fact is important for understanding the smooth rigidity phenomena, studied in the present paper, so in Section \ref{Sec:an.example} below we present in some detail an example.

\medskip

It is easy to see that there is {\it a uniform upper bound for the $d$-rigidity of all the subsets $Z\subset \hat B^n$.} Indeed, consider a certain $C^\infty$ function $\phi$, which vanishes identically on $\hat B^n$ and satisfies $M_0(\phi)=1$. Then $\psi$ vanishes in $Z$, and hence ${\cal RG}_d(Z)\le M_{d+1}(\phi).$

\medskip

Another simple observation is the following:

\bp\label{prop:Z.interior}
For any $Z \subset \hat B^n$ with a non-empty interior,
$$
{\cal RG}_d(Z) \ge \frac{(d+1)!}{2^{d+1}}.
$$
\ep
This fact easily follows from Proposition \ref{prop:d.points}. We just restrict any function $f\in U_d(Z)$ to a certain straight line $\ell$, passing through $z_0$ with $|f(z_0)|=1,$ and through an interior point of $Z$. The present paper extend this result to all sufficiently dense $Z$.

\medskip

Let us mention also an old result of \cite{Yom1}, related to smooth rigidity. Informally it can be stated as follows: if the set of zeros $Y(f)$ of a smooth function $f$ on $B^n$ does not look like a union of smooth hypersutfaces, of a total area bounded by a constant, depending only on $n,d$, then the norm of $f^{(d+1)}$ is not smaller than a certain positive constant, depending only on $n,d$.

\medskip

Many of the ``near-polynomiality'' results of \cite{Yom2,Yom.Com} can be naturally interpreted in terms of smooth rigidity. We plan to present some new results in this direction separately.

\subsection{Remez constant of $Z$}\label{Sec:Remez.zero.sets}

Another ingredient we need is a definition and some properties of the Remez (or Lebesgue, or norming, ...) constant (see, e.g. \cite{Rem,Bru.Yom,Yom3} and references therein).

\bd\label{Remez.constant}
For a set $Z\subset B^n \subset {\mathbb R}^n$ the Remez constant ${\cal R}_d(Z)$ is the minimal $K$
for which the inequality
$$
\sup_{B^n}\vert P \vert \leq K \sup_{Z}\vert P \vert
$$
is valid for any real polynomial $P(x)=P(x_1,\dots,x_n)$ of degree $d$.
\ed
Clearly, we always have ${\cal R}_d(Z)\ge 1.$ For some $Z$ the Remez constant ${\cal R}_d(Z)$ may be equal to $\infty$. In fact, ${\cal R}_d(Z)$ is infinite if and only if $Z$ is contained in the set of zeroes
$$
Y_P=\{x\in {\mathbb R}^n, \ | \ P(x)=0\}
$$
of a certain polynomial $P$ of degree $d$. Sometimes it is convenient to use the inverse Remez constant $\hat {\cal R}_d(Z):=\frac{1}{{\cal R}_d(Z)}.$

\subsection{Remez constant and rigidity (\cite{Yom5})}\label{Sec:main.results}

In \cite{Yom4,Yom5} we show that the rigidity ${\cal RG}_d(Z)$ and the Remez constant $\hat {\cal R}_d(Z)$ are closely connected:

\bt\label{thm:main1}(\cite{Yom5})
For any $Z \subset \hat B^n$, \ \ \ $\frac{(d+1)!}{2}\hat {\cal R}_d(Z)\le {\cal RG}_d(Z)$.
\et
This lower bound is not always sharp. Indeed, by Proposition \ref{prop:Z.interior}, for any $Z \subset \hat B^n$ with a non-empty interior, we have ${\cal RG}_d(Z) \ge \frac{(d+1)!}{2^{d+1}},$ while $\hat {\cal R}_d(Z)$ can be arbitrarily small. However, the bound of Theorem \ref{thm:main1} is sharp, up to constants, for finite sets with a controlled minimal distance between the points:

\bt\label{thm:main2}(\cite{Yom5})
Let $Z \subset B^n$ be a finite set, and let $\rho$ be the minimal distance between the points of $Z$. Then
$$
\frac{1}{2}\hat {\cal R}_d(Z)\le {\cal RG}_d(Z)\le \frac{C(n,d)}{\rho^{d+1}}\hat {\cal R}_d(Z).
$$
\et
This theorem can be considered as one of possible generalizations of Proposition \ref{prop:d.points} to higher dimensions.

\medskip

We use in \cite{Yom5} a Remez-type inequality for discrete sets, obtained in \cite{Yom3}, in combination with Theorem \ref{thm:main1}, to provide the following lower bound for the rigidity in terms of the density of $Z$:



\bt\label{thm:main3}(\cite{Yom5})
Let $Z \subset B^n$ be a finite set, and let $\rho$ be the minimal distance between the points of $Z$. Assume that the cardinality $|Z|$ satisfies
$|Z| > (4d)^n(\frac{1}{\rho})^{n-1}$. Then
$$
0 < \frac{(d+1)!}{2}\left ( \frac{|Z|\rho^n- (4d)^n\rho}{4n}\right )^d \le {\cal RG}_d(Z).
$$
\et
Both Theorem \ref{thm:main3} and the results of the present paper provide lower bounds for the smooth rigidity ${\cal RG}_d(Z)$ in terms of certain ``densities'' of $Z$, expressed via the covering number $M(Z,\e)$. 

\medskip

To stress the difference of the density-based lower bounds on ${\cal RG}_d(Z)$ via the Remez constant, as in \cite{Yom5}, and the results of the present paper, we consider the following example: let $Z$ be a finite $h$-dense sets in the $s$-cube $Q^n_s, \ \ 0<h<s\le \frac{1}{10}$. Assume also that the minimal distance between the points of $Z$ is $\rho=\frac{h}{3}$.

\bp\label{prop:h.dense.Remez}
If $h\le K_1s^{n}$, with $K_1=2\cdot 3^{n-1}(4d)^n$, then
$$
{\cal RG}_d(Z) \ge \frac{(d+1)!}{2}(\frac{s}{12n})^{nd}.
$$
\ep
\pr
Since $Z$ is $h$-dense, the cardinality $|Z|$ is at least $(1+o(h))(\frac{s}{h})^n.$ Therefore the condition $h\le K_1s^n$ implies
$$
|Z| \ \ge \ (1+o(h))(\frac{s}{h})^n = (1+o(h))\frac{s^n}{h}\cdot (\frac{1}{h})^{n-1}\ge (1+o(h))K_1(\frac{1}{h})^{n-1}= 
$$
$$
= (1+o(h))2(4d)^n(\frac{3}{h})^{n-1} \ \ = \ \ (1+o(h))2(4d)^n(\frac{1}{\rho})^{n-1}.
$$
We conclude that
$$
|Z|\rho^n \ge (1+o(h))2(4d)^n\rho,
$$
and hence
$$
|Z|\rho^n- (4d)^n\rho \ \ge \ \frac{1}{3}|Z|\rho^n \ \ge \ \frac{1}{3}(\frac{s}{h})^n(\frac{h}{3})^n \ = \ (\frac{s}{3})^n.
$$
Substituting this bound into the expression of Theorem \ref{thm:main3}, we finally get
$$
{\cal RG}_d(Z)\ge \frac{(d+1)!}{2}(\frac{s}{12n})^{nd}.
$$
This completes the proof of Proposition \ref{prop:h.dense.Remez}.


\medskip

Now, in Theorem \ref{main:rigid.h.regular} we require a rather high density of $Z$, namely, $h\le Ks^{n(d+1)}$. But the conclusion is that ${\cal RG}_d(Z)\ge d!$, independently of the size $s$ of the cube $Q^n_s$, containing $Z$.

\smallskip

In Proposition \ref{prop:h.dense.Remez} the density requirement is much weaker: $h\le K_1s^{n}$, and the conclusion is also weaker: ${\cal RG}_d(Z)\ge \frac{(d+1)!}{2}(\frac{s}{12n})^{nd}$. In particular, as the size $s$ of the cube tends to zero, this bound decreases as $s^{nd}$.

\subsection{A set $Z$ with a small positive $d$-rigidity}\label{Sec:an.example}

As it was mentioned above, in dimensions higher than one there is no jump in the possible values of $d$-rigidity. For $Z\subset \hat B^n, \ n\ge 2,$ the $d$-rigidity ${\cal RG}_d(Z)$ may attain arbitrarily small positive values (\cite{Yom5}).

\medskip

Extending \cite{Yom5} we give here an explicit example for $n=2, \ d=1$ . Consider a plane triangle $Z_h$, defined as
$$
Z_h=\{(-\frac{1}{2},0),(0,h),(\frac{1}{2},0)\}, \ 0<h\le \frac{1}{8}.
$$
We'll show that $\frac{h}{4}\le {\cal RG}_1(Z_h)\le 4h.$

\smallskip

Let us first estimate the Remez constant ${\cal R}_1(Z_h)$. Assume that a first degree polynomial $P(x,y)=ax+by+c$ satisfies $|P|\le 1$ on $Z_h$. Then, in particular, $|P(-\frac{1}{2},0)|=|-\frac{1}{2}a+c|\le 1$, \ $|P(\frac{1}{2},0)|=|\frac{1}{2}a+c|\le 1$. We conclude that $|P(-\frac{1}{2},0)+P(\frac{1}{2},0)|=|2c|\le 2$, or $|c|\le 1$. In the same way we get $|P(-\frac{1}{2},0)-P(\frac{1}{2},0)|=|a|\le 2$.

\medskip

Next we have $|P(0,h)|=|bh+c|\le 1$, or $|bh|\le |c|+1\le 2$, or $|c|\le \frac{2}{h}$. Therefore, on the square $|x|,|y|\le 1$, which contains the disk $B^2$, we obtain $|P(x,y)|\le 1+2+\frac{2}{h}=\frac{3h+2}{h}<\frac{4}{h}$. Thus
$${\cal R}_1(Z_h)\le\frac{4}{h}, \ \ \hat {\cal R}_1(Z_h)\ge\frac{h}{4}.
$$
By Theorem \ref{thm:main1} we get
$$
{\cal RG}_1(Z_h)\ge \frac{(1+1)!}{2}\frac{h}{4}=\frac{h}{4}.
$$

In order to obtain an upper bound for ${\cal RG}_1(Z_h)$ we construct a polynomial $Q(x,y)$ of degree two, vanishing on $Z_h$, and estimate its maximum on $B^2$ and its second derivative.

\medskip

Put $Q(x,y)=y+4h(x^2-\frac{1}{4})$. Then $Q$ vanishes on $Z_h$, the norm of its second derivative is $2h$, and $1-4h \le M_0(Q) \le 1+4h$. After we normalize $Q$ to $\tilde Q$ with $M_0(\tilde Q)=1$, we still have $M_2(\tilde Q)\le \frac{2h}{1-4h} \le 4h$. We conclude that
$$
\frac{h}{4}\le {\cal RG}_1(Z_h)\le 4h.
$$
In particular, for $h\to 0$ the $1$-rigidity of $Z_h$ tends to $0$, remaining positive.

\subsection{Smooth rigidity via topology}\label{Sec:topol.rigidity}

Let us mention that in recent papers \cite{Yom6} and (in other terms) in \cite{Ler.Ste} rigidity inequalities are obtained based on topology of $Z$. The approach of the present paper is different, but not independent of the topological approach. Indeed, for the zero set $Z$ of $f$ consisting of $d+1$ nested hypersurfaces, we easily find straight lines $\ell$ crossing $Z$ at $d+1$ points. Are there topological (or combined topological-geometric) conditions on $Z$, providing curves with small high-order derivatives, crossing $Z$ at $d+1$ points?

\medskip

From now on we turn to the proof of Theorems \ref{thm:main.intro1} and \ref{main:rigid.h.regular}.


\section{Comparing derivatives of $f(x)$ and $f(\o(t))$}\label{Sec:Comparison.Deriv}
\setcounter{equation}{0}

In this section we prove an inequality, which compares the highest derivatives of a function $f$ on $B^n$, and of its restriction to a given smooth curve $\o$. More accurately, we consider $C^{d+1}$-smooth curves $\o:[-1,1]\to B^n$, given in the coordinate form by
$$
\o(t)=(\o_1(t),\ldots, \o_n(t)),
$$
and compare the derivatives of $f$ and of the composition $f\circ \o$. The ``chain rule'' expressions for higher derivatives are rather complicated, so in the next section we recall these expressions, and summarise the required facts about them.

\subsection{Derivatives of $f(\o(t))$: symbolic expressions}\label{Sec:Composition.Deriv}

We consider the composition $g(t)=f(\o(t))$. Thus, $g(t)$ is a $C^{d+1}$-smooth function of one real variable $t\in [-1,1]$. By the formula for the derivatives of the composition we have:

$$
g'(t)=f'(\o(t))\o'(t),
$$
$$
g''(t)=f''(\o(t))\o'(t)^2+f'(\o(t))\o''(t),
$$
$$
g'''(t)=f'''(\o(t))\o'(t)^3+3f''(\o(t))\o'(t)\o''(t)+f'(\o(t))\o'''(t).
$$
We continue this list till the fifth derivative, omitting in the further formulas the arguments $t$ and $\o(t)$:
$$
g^{(iv)}=f^{(iv)}\o'^4+6f'''\o'^2\o''+f''(3\o''^2+4\o'\o''')+f'\o^{(iv)},
$$
$$
g^{(v)}=f^{(v)}\o'^5+ 10f^{(iv)}\o'^3\o''+f'''(15\o'\o''^2+5\o'^2\o''') + f''(10\o''\o'''+5\o'\o^{(iv)})+f'\o^{(v)}.
$$
$$
..................
$$
For $n=1$ the formulas above do not need further interpretation. In higher dimensions $n$, we obtain an expression which must be interpreted in terms of tensor calculus of high order derivatives:

\be\label{eq:composition}
g^{(d+1)}(t)= \sum_{k=1}^{d+1}f^{(k)}(\o(t))\ast\kappa_k(\o(t)),
\ee
where $\kappa_{d+1}(\o(t))=(\o'(t))^{d+1}$, $\kappa_{1}(\o(t))=\o^{(d+1)}(t)$, and $\kappa_{k}(\o(t))$ are polynomials in the derivatives of $\o(t)$. These polynomials become homogeneous, of degree $d+1$, if we assign the variable $\o^{(k)}(t)$ the degree $k$. The star product $\ast$ denotes the appropriate tensor product.

\smallskip

For an explicit form in higher dimensions, consider, for example, the case $n=2, \ d=1,2$. Here $\o(t)$ is given by two functions $\o(t)=(x(t),y(t))$, and for $g(t)=f(\o(t))$ we have
$$
\frac{dg(t)}{dt}=\frac{\partial f(\o(t))}{\partial x}\frac{dx}{dt}+\frac{\partial f(\o(t))}{\partial y}\frac{dy}{dt},
$$
$$
\frac{d^2g(t)}{dt^2}=\frac{\partial^2 f(\o(t))}{\partial x^2}(\frac{dx}{dt})^2+2\frac{\partial f(\o(t))}{\partial x\partial y}\frac{dx}{dt}\frac{dy}{dt}+\frac{\partial^2 f(\o(t))}{\partial y^2}(\frac{dy}{dt})^2+
$$
$$
+\frac{\partial f(\o(t))}{\partial x}\frac{d^2x}{dt^2}+\frac{\partial f(\o(t))}{\partial y}\frac{d^2y}{dt^2},
$$
et cetera. For our purposes it is enough to state the following properties of the expressions above:

\medskip

\noindent 1. $g^{(d+1)}(t)$ is a polynomial in the partial derivatives of $f$ and in the derivatives of $\o$. This polynomial is linear in $f$, and it becomes homogeneous, of degree $d+1$, in the derivatives of $\o$, if we assign the variables $\o_i^{(k)}(t)$ the degree $k$. The number of the monomials in this polynomial does not exceed $B_1(d,n)$ and the coefficients with the monomials do not exceed $B_2(d,n)$.

\medskip

\noindent 2. The only term in the expression for $g^{(d+1)}(t),$ given by (\ref{eq:composition}), which contains the highest ($d+1$-st) order derivatives of $f$, is $f^{(d+1)}(\o(t))\ast\kappa_{d+1}(\o(t))$. This is a homogeneous form $Q$ of degree $d+1$ in the first order derivatives of $\o$, with the coefficients - all the partial derivatives of $f$ of order $d+1$, i.e.
\be\label{eq:high.homog}
f^{(d+1)}(\o(t))\ast\kappa_{d+1}(\o(t))=Q(d^{d+1}f,\frac{d\o}{dt})=\sum_{|\alpha|=d+1}\frac{\partial^\alpha f(\o(t))}{\partial x^\alpha}(\frac{d\o}{dt})^\alpha,
\ee
where $\frac{d\o}{dt}=(\frac{dx_1}{dt},\ldots,\frac{dx_n}{dt}),$ in the usual multi-index notations.

\medskip

\noindent 3. Each monomial in the polynomials $\kappa_k(\o(t)), \ k=1,\ldots,d,$ in (\ref{eq:composition}) necessarily contains derivatives of $\o$ of order at least two.


\subsection{Derivatives of $g(t)=f(\o(t))$: main inequality}\label{Sec:Composition.Deriv1}

Let $f$ be a $C^{d+1}$-smooth function on $B^n$. As above, we use the norm $||d^kf(z)||$ of the derivatives of functions $f$ on $B^n$, defined as
$$
||d^kf(z)||=\max_{|\alpha|=k} \ |\frac{\partial^\alpha f(z)}{\partial x^\alpha}|.
$$
It will be important to consider separately the highest derivative $f^{(d+1)}$ and the derivatives $f^{(k)}, \ k=1,\ldots,d.$ In particular, for each $z\in B^n$ we put
$$
\mu_d(f,z):= \max_{k=1}^d ||d^kf(z)||,  \  \mu_d(f)=\max_{z\in B^n} \mu_d(f,z).
$$
Now we consider $C^{d+1}$-smooth curves $\o:[-1,1]\to B^n$, given in the coordinate form by
$$
\o(t)=(\o_1(t),\ldots, \o_n(t)).
$$
Also here we use the max-norm $||d^k\o(t)||$ of the derivatives of $\o$, defined as
$$
||d^k\o(t)||=\max_{i=1}^n \ |\o_i^{(k)}(t)|.
$$
For curves $\o$ it will be important to consider separately the velocity $\o'(t)$, and the higher order derivatives of $\o$. This is because we want to measure the deviation of $\o$ from a straight line. In particular, for each $t\in [-1,1]$ we put
$$
\nu_d(\o,t)=\max_{k=2}^{d+1}||d^k\o(t))||, \  \  \nu_d(\o)=\max_{t\in [-1,1]}\nu_d(\o,t).
$$


\bp\label{prop:composition}
Let $f$ be a $C^{d+1}$-smooth function on $B^n$, and let $\o$ be a $C^{d+1}$-smooth curve in $B^n$, satisfying
$$
||d\o(t)|| \le 1, \ \nu_d(\o) \le 1.
$$
Then for $g(t)=f(\o(t))$ and for each $t\in [-1,1]$ we have
$$
||d^{d+1}f(\o(t))||\ge \ \frac{C_1(n,d)}{||\o'(t)||^{d+1}}\left (|g^{(d+1)}(t)|-C_2(d,n)\mu_d(f,\o(t))\nu_d(\o,t)\right ),
$$
with $C_1(d,n),C_2(n,d)$ depending only on $d$ and $n$.
\ep
\pr
From (\ref{eq:composition}) we get
\be\label{eq:composition1}
f^{(d+1)}(\o(t))\ast\kappa_{d+1}(\o(t))=g^{(d+1)}(t)-\sum_{k=1}^{d}f^{(k)}(\o(t))\ast\kappa_k(\o(t)).
\ee
Noice that this is a scalar expression. Put $\Sigma:=\sum_{k=1}^{d}f^{(k)}(\o(t))\ast\kappa_k(\o(t)).$

\bl\label{lem:lower.ord.terms}
Under the assumptions of Proposition \ref{prop:composition} we have
$$
|\Sigma|\le C_2(d,n)\mu_d(f,\o(t))\nu_d(\o,t),
$$
with $C_2(d,n)$ being a constant depending only on $n$ and $d$.
\el
\pr
By the definitions, the norms of the partial derivatives of $f$ at $\o(t)$ satisfy $||f^{(k)}(\o(t))||\le \mu_d(f,\o(t)), \ k=1,\ldots,d.$ Next,  by the property (3) above, each monomial of $\kappa_k(\o(t)), \ k=1,\ldots,d,$ contains a certain derivative of $\o$ of order at least two, which is bounded in absolute value by $\nu_d(\o,t)$. The rest of the terms in the monomial are bounded in absolute value either by $1$, or by $\nu_d(\o)\le 1$, by the assumptions. We conclude that each monomial is bounded in absolute value by $\mu_d(f,\o(t))\nu_d(\o,t)$. Finally, by (1), the number of the monomials in $\kappa_k(\o(t))$ does not exceed $B_1(d,n)$ and the coefficients with the monomials do not exceed $B_2(d,n)$. The tensor product $\ast$ ads at most a coefficient $B_3(n,d)$. This  completes the proof of Lemma \ref{lem:lower.ord.terms}, with  $C_2(d,n)=d B_1(n,d)B_2(d,n)B_3(d,n).$    $\square$.

\medskip

Now we use expression (\ref{eq:high.homog}) for the left hand side of (\ref{eq:composition1}).

\bl\label{lem:leading.term}
\be\label{eq:composition2}
|f^{(d+1)}(\o(t))\ast\kappa_{d+1}(\o(t))| \le B_4(n,d)||f^{(d+1)}(\o(t)||\cdot ||\o'(t)||^{d+1}.
\ee
\el
\pr
By our definitions of the norms, the partial derivatives of order $d+1$ of $f$ at $\o(t)$ do not exceed $||f^{(d+1)}(\o(t)||$, while the first derivatives of $\o(t)$ do not exceed $||\o'(t)||$. The result now follows directly from (\ref{eq:high.homog}). $\square$

\smallskip

Hence we obtain, via Lemmas \ref{lem:lower.ord.terms} and \ref{lem:leading.term},
$$
B_4(n,d)||f^{(d+1)}(\o(t))||\cdot ||\o'(t)||^{d+1}\ge |f^{(d+1)}(\o(t))\ast\kappa_{d+1}(\o(t))| =|g^{(d+1)}(t)-\Sigma|\ge
$$
$$
\ge |g^{(d+1)}(t)|-|\Sigma|\ge |g^{(d+1)}(t)|-C_2(d,n)\mu_d(f,\o(t))\nu_d(\o,t).
$$
Dividing by $B_4(n,d)||\o'(t)||^{d+1}$ we finally get the required inequality
$$
||d^{d+1}f(\o(t))||\ge \ \frac{C_1(n,d)}{||\o'(t)||^{d+1}} ( |g^{(d+1)}(t)|-C_2(d,n)\mu_d(f,\o(t))\nu_d(\o,t) ),
$$
with $C_1(d,n)=\frac{1}{B_4(n,d)}$. This completes the proof of Proposition \ref{prop:composition}. $\square$

\bc\label{cor:compar.compos}
For $f$ and $\o$ as above we have
\be\label{eq:composition3}
M_{d+1}(f)\ge C_1 \left ( M_{d+1}(g) - C_2 \mu_d(f)\nu_d(\o) ) \right.
\ee
\ec
\pr
We fix $t_0$ to be one of the points, where the maximum $M_{d+1}(g)$ is attained, and apply Proposition \ref{prop:composition} at $t_0$, taking into account that by the assumptions $||\o'||\le 1$. $\square$

\smallskip

Next we want to exclude the dependence of the right-hand side of (\ref{eq:composition3}) on $\mu_d(f)$. We can do it for $\nu_d(g)$ sufficiently small, i.e. for curves $\o$ with sufficiently small high order derivatives. Let us recall a simple bound for ``intermediate derivatives'' (see, e.g. \cite{Yom6}, Lemma 7.1):

\bl\label{lem:norms.der}
Let $f$ be a $C^{d+1}$-smooth function on $B^n.$ Then for $k=1,2,\ldots,d$ we have
$$
M_k(f)\le B_5(n,d)M_0(f)+ B_6(n,d)M_{d+1}(f).
$$
\el

\bp\label{prop:composition6}
For $f$ and $\o$ as above assume that $M_0(f)=1$ and assume that $\nu_d(\o)\le C_3(n,d)$, where the constant $C_3$ is defined below. Then we have
\be\label{eq:composition4}
M_{d+1}(f)\ge \frac{1}{2}C_1 (M_{d+1}(g) - \frac{1}{10}).
\ee
\ep
\pr
Applying Lemma \ref{lem:norms.der} to $f$ with $M_0(f)=1$, we obtain
$$
M_k(f)\le B_5+ B_6M_{d+1}(f), \  k=1,2,\ldots,d,
$$
and hence $\mu_d(f)\le B_5+ B_6M_{d+1}(f)$. From Corollary \ref{cor:compar.compos} we conclude that
$$
M_{d+1}(f)\ge C_1(M_{d+1}(g) - C_2 (B_5+ B_6M_{d+1}(f))\nu_d(\o)) =
$$
$$
=C_1M_{d+1}(g) - (B_7 + B_8M_{d+1}(f))\nu_d(\o),
$$
where $B_7=C_1C_2B_5, \ B_8=C_1C_2B_6$. We rewrite this expression as
$$
(1+B_8\nu_d(\o))M_{d+1}(f) \ge C_1M_{d+1}(g) - B_7 \nu_d(\o).
$$
Now we assume that $\nu_d(\o)\le C_3(n,d):=\min \{\frac{1}{2B_7}, \frac{1}{2B_8}, \frac{C_1}{10B_7} \}$. Under this assumption we get
$$
2M_{d+1}(f) \ge C_1(M_{d+1}(g) - \frac{1}{10}),
$$
or
$$
M_{d+1}(f) \ge \frac{1}{2}C_1(M_{d+1}(g) - \frac{1}{10}).
$$
This completes the proof of Proposition \ref{prop:composition6}. $\square$

\subsection{Some conclusions for $d$-rigidity}\label{sec:concl.d.rigid}

\bt\label{thm:composition.rigid}
Let $f$ be a $C^{d+1}$-smooth function on $B^n$, with $M_0(f)=1$, and let $\o:[-1,1]\to B^n$ be a $C^{d+1}$-smooth curve in $B^n$, satisfying
$$
||d\o(t)|| \le 1, \ t\in [-1,1], \ \nu_d(\o) \le C_3.
$$
Assume that the curve $\o$ passes through a certain point $z_0\in S^n$ with $|f(z_0)|=1$, and through some points $z_1,\ldots,z_{d+1}$ with $f(z_j)=0, \ j=1,\ldots,d+1$. Then
$$
M_{d+1}(f) \ge C_4(n,d):=\frac{1}{2}C_1(n,d)(\frac{(d+1)!}{2^{d+1}} - \frac{1}{10})>0.
$$
\et
\pr First we apply Proposition \ref{prop:d.points} to $g(t)=f(\o(t))$, and conclude that $M_{d+1}(g)\ge \frac{(d+1)!}{2^{d+1}}$. Next we apply to $f$ and $\o$ Proposition \ref{prop:composition6}. This completes the proof. $\square$

\medskip

In what follows we look for geometric conditions on the zero set $Z$ of $f$ which imply existence of the curves $\o$, satisfying conditions of Theorem \ref{thm:composition.rigid}.

\section{$d$-thickness of zero sets}\label{Sec:thickness}
\setcounter{equation}{0}

Our next goal is to introduce a certain geometric characteristic of closed subsets $Z\subset \hat B^n$, which we call a $d$-thickness of $Z$. It estimates the ``size'' of $Z$ with respect to a possibility to draw a smooth curve $\o$ with small high-order derivatives, through a given point $z_0\in S^n$ and some $d+1$ distinct points of $Z$.

\subsection{Definition of $d$-thickness}\label{sec:curve.dens}

Let us recall that for a $C^{d+1}$-smooth curve $\o:[-1,1]\to B^n$ we put
$$
\nu_d(\o,t)=\max_{k=2}^{d+1}||d^k\o(t))||, \  \  \nu_d(\o)=\max_{t\in [-1,1]}\nu_d(\o,t).
$$

\bd\label{def:cutting.lines}
Let a set $Z\subset \hat B^n$, and a point $z_0\in S^n$ be given. We denote $\O_d(Z,z_0)$ the collection of all $C^{d+1}$-smooth curves $\o:[-1,1]\to B^n,$ passing through $z_0$, and through certain $d+1$ distinct points of $Z$, which satisfy the following conditions:

\smallskip

\noindent 1. The velocity $\frac{d\o(t)}{dt}$ of $\o$ satisfies $1\ge ||\frac{d\o(t)}{dt}|| > 0, \ t\in [-1,1].$

\smallskip

\noindent 2. $\nu_d(\o)\le 1.$

\medskip

The $d$-thickness $\nu_d(Z,z_0)$ is the minimum over all $\o\in \O_d(Z,z_0)$ of $\nu_d(\o)$.

\medskip

The $d$-thickness $\nu_d(Z)$ is the maximum over all $z_0\in S^n$ of $\nu_d(Z,z_0)$.
\ed

If we can find a straight line $\ell$ in ${\mathbb R}^n,$ passing through $z_0$, and through certain $d+1$ distinct points of $Z$, then $\nu_d(Z)=0$. Otherwise, $\nu_d(Z)$ measures the minimal ``high order curvature'' (taking into account high order derivatives) of the smooth curves, passing through $z_0$, and through certain $d+1$ distinct points of $Z$.

\medskip

A direct consequence of Definition \ref{def:cutting.lines} and of Theorem \ref{thm:composition.rigid} is the following result:

\bt\label{thm:rigidity.width}
For $Z\subset \hat B^n$ if we have $\nu_d(Z)\le C_3(n,d)$, then
$$
{\cal RG}_d(Z) \ge C_4(n,d)>0.
$$
\et
\pr
Let us recall that for $Z\subset \hat B^n$ we denote $U_d(Z)$ the set of all $C^{d+1}$ smooth functions $f(z)$ on $B^n$, vanishing on $Z$, with $M_0(f)=1$.

\medskip

Let now a function $f\in U_d(Z)$ be given, and let $z_0$ be a point where $|f(z_0)|=1$. Since $\nu_d(Z)\le C_3(n,d)$, we can find a curve $\o$ passing through $z_0\in S^n$ and through some points $z_1,\ldots,z_{d+1}$ in $Z$, and satisfying
$$
||d\o(t)|| \le 1, \ t\in [-1,1], \ \nu_d(\o) \le C_3.
$$
By Theorem \ref{thm:composition.rigid} we conclude that $M_{d+1}(f) \ge C_4(n,d)$. Then via Definition \ref{def:rigidity} we obtain
$$
{\cal RG}_d(Z)=\inf_{f\in U_d(Z)}M_{d+1}(f)\ge C_4(n,d).
$$
This completes the proof of Theorem \ref{thm:rigidity.width}. $\square$

\medskip

In particular, for sets $Z$ with a small $d$-thickness $(\nu_d(Z)\le C_3(n,d))$, like for sets with a non-empty interior, the $d$-rigidity ${\cal RG}_d(Z)$ is uniformly bounded from below by $C_4(n,d),$ independently of the size of $Z$.

\subsection{An approximate $ad$-thickness}\label{sec:lin.dens1}

In what follows, we estimate the $d$-thickness $\nu_d(Z)$, essentially, via bounding another, somewhat simpler geometric characteristic of sets $Z$, which we call an approximate $ad$-thickness $\bar \nu_d(Z)$. The main advantage of the approximate $ad$-thickness is that it does not involve directly high-order derivatives, and involves only straight lines and distances. However, in the present paper we do not try to further develop this research direction, planning to present more results separately.

\medskip

For each straight line $\ell$ in ${\mathbb R}^n$, and for each collection ${\cal Z}$ of $d+1$ distinct points $z_1,\ldots,z_{d+1}$ in ${\mathbb R}^n$ we put $\rho(\ell,{\cal Z})$ to be the maximal distance of the points in ${\cal Z})$ to $\ell$.

\smallskip

We also put $\kappa(\ell,{\cal Z})$ to be the minimal distance between the projections of the points in ${\cal Z})$ to $\ell$.

\smallskip

Now we define $\mu_d(\ell,{\cal Z})$ as

\be\label{def.ad.th}
\mu_d(\ell,{\cal Z})=\frac{\rho(\ell,{\cal Z})}{\kappa(\ell,{\cal Z})^d}.
\ee
Now let a point $z_0\in S^n$ be given. Consider the set ${\cal L}(z_0)$ of all the straight lines $\ell$ in ${\mathbb R}^n$ passing through $z_0$ and intersecting the ball $\hat B^n$.

\bd\label{def:cutting.lines1}
Let a set $Z\subset \hat B^n$, and a point $z_0\in S^n$ be given. The point-wise approximate $ad$-thickness $\bar \nu_d(Z,z_0)$ is defined as
$$
\bar \nu_d(Z,z_0)= \min_{\ell\in {\cal L}(z_0)} \ \min_{{\cal Z}\subset Z} \ \mu_d(\ell,{\cal Z}).
$$
The approximate $ad$-thickness $\bar \nu_d(Z)$ is defined as
$$
\bar \nu_d(Z)=\max_{z_0\in S^n} \ \bar \nu_d(Z,z_0).
$$
\ed
Our goal now is to prove that the approximate $ad$-thickness $\bar \nu_d(Z)$ bounds from above the $d$-thickness $\nu_d(Z)$.

\subsection{$\bar \nu_d(Z)$ bounds $\nu_d(Z)$}\label{sec:lin.dens2}

To get this bound (under certain restrictions) it is enough to prove the following result:

\bp\label{prop:compar.nu.bar.nu}
Let a point $z_0\in S^n$, a straight line $\ell\in {\cal L}(z_0)$, and a collection ${\cal Z}$ of $d+1$ distinct points $z_1,\ldots,z_{d+1}$ in $\hat B^n$ be given. Assume that $\rho=\rho(\ell,{\cal Z})$ and $\kappa=\kappa(\ell,{\cal Z})$ satisfy $\frac{1}{10} > \kappa \ge 10\rho$.

\medskip

Then there exists a smooth curve $\o:[-1,1]\to B^n$, passing through $z_0$ and through $z_1,\ldots,z_{d+1}$, and satisfying
$$
||d\o(t)|| \le 1, \ t\in [-1,1], \ \nu_d(\o) \le \bar D(d) \cdot \mu_d(\ell,{\cal Z}),
$$
with the constant $\bar D(d)$ depending only on $d$.
\ep
The proof of this proposition is given in the next section. As a direct consequence we obtain:

\bc\label{cor:compar.nu.bar.nu}
For each $Z\subset \hat B^n$ we have
$$
\nu_d(Z)\le \bar D(d) \bar \nu_d(Z).
$$
\ec

Our geometric calculations in Section \ref{sec:entr.bounds.nu} below will be devoted, essentially, to bounding $\nu_d(Z)$ in terms of covering density of $Z$.

\subsection{Smooth curves through given points}\label{Sec:curves.through}

In fact, we give here a more detailed statement of the required result. Proposition \ref{prop:straight.curved} below covers a situation, where a straight line $\ell$ passes close enough to certain $d+1$ points $z_1,\ldots,z_{d+1}$. We deform $\ell$ into a curve with controlled high-order derivatives, which passes exactly through $z_1,\ldots,z_{d+1}$.

\bp\label{prop:straight.curved}
Let ${\cal Z}=\{z_1,\ldots,z_{d+1}\}\subset \hat B^n$, and $z_0\in S^{n-1},$ be given. Assume that there exists a straight line $\ell,$ passing through $z_0$, such that the following conditions are satisfied:

\medskip

\noindent 1. The distance of each point $z_i$ to $\ell, \ i=1,\ldots,d,$ is at most $\rho>0$.

\medskip

\noindent 2. The distance between the projections $\eta_i$ of the points $z_i$ to $\ell$ is at least $\kappa>0$. As above, we assume that        $\frac{1}{10} > \kappa \ge 10\rho$.

\medskip

Then there exists a $C^{d+1}$-smooth curve $\o: [-1,1]\to B^n,$ passing through the points $z_0,z_1,\ldots,z_{d+1}$ and satisfying
$$
M_1(\o) \ \le \ 1, \ \ \ M_k(\o) \ \le \ \frac{D_k\rho}{\kappa^{k}}, \ \ \ k=2,3,\ldots, d+1,
$$
with $D_k$ being the constants, depending only on $k$. In particular,
$$
\nu_d(\o)\le \bar D(d) \frac{\rho}{\kappa^d} = \bar D(d) \cdot \mu_d(\ell,{\cal Z}),
$$
with $\bar D(d)=\max_{k=2,\ldots,d+1} \ D_k$.
\ep
\pr
We can assume that the projections of the points $z_i$ to $\ell$ are ordered. Denote by $v$ the unit vector in the direction of $\ell$, pointing towards the projections of ${\cal Z}$. We denote by $v_i, \ i=1,\ldots,d+1,$ the vectors, orthogonal to $\ell$, connecting the projections $\eta_i$ of the points $z_i$ to $\ell$, with the points $z_i$ themselves.

\smallskip

First we consider the Euclidean length $t$ as the coordinate on $\ell$, with the origin at the point $z_0$. Let $\tau_i$ be the $t$-coordinates of the projections $\eta_i$.

\medskip


\medskip

Now the curve $\tilde \o$ is defined as
$$
\tilde \o(t)=z_0+tv + \sum_{i=1}^{d+1}v_i\phi(\frac{2(t-\tau_j)}{\kappa}),
$$
where $\phi(t)$ is a fixed $C^\infty$ ``Gaussian-like'' function on $[-1,1]$, identically equal to zero near the ends of $[-1,1]$, and outside of this interval, with $\phi(0)=1$. Since by the assumptions the distance between $\tau_i$ is at least $\kappa$, the supports of the functions $\phi(\frac{2(t-\tau_i)}{\kappa})$ are disjoint, and thus we have $\o(\tau_i)=z_i, \ i=0,1,\ldots,d+1.$

\medskip

For the derivatives of $\tilde \o(t)$ we obtain
$$
\frac{d\tilde \o(t)}{dt}=v+\sum_{i=1}^{d+1}v_i\phi'(\frac{2(t-\tau_i)}{\kappa})\cdot \frac{2}{\kappa},
$$
$$
\frac{d^k\tilde \o(t)}{dt^k}=\sum_{i=1}^{d+1}v_i\phi^{(k)}(\frac{2(t-\tau_i)}{\kappa})\cdot (\frac{2}{\kappa})^k, \ \ \ k\ge 2.
$$
Once more, since by the assumptions the distance between $\tau_i$ is at least $\kappa$, the supports of the functions $\phi(\frac{2(t-\tau_i)}{\kappa})$ are disjoint, and we obtain
$$
M_1(\tilde \o)\le 1+\max_{j=1}^{d+1}||v_i||M_1(\phi)\cdot \frac{2}{\kappa}\le 1+2M_1(\phi)\cdot \frac{\rho}{\kappa},
$$
$$
M_k(\tilde \o)\le 2^k M_k(\phi)\cdot \frac{\rho}{\kappa^k}, \ \ k\ge 2.
$$
Now we perform an affine reparametrization of the curve $\o$, mapping the interval $[-1,1]$ with the coordinate $\eta$ via an affine mapping $t=\psi(\eta)$ onto the interval $[-\frac{1}{10}, \tau_{d+1}+\frac{1}{10}]$ in the coordinate $t$ on $\ell$.

\smallskip

Finally we put $\o(\eta)=\tilde \o(\psi(\eta), \ \eta \in [-1,1]$. Since the length of the interval $[-\frac{1}{10}, \tau_{d+1}+\frac{1}{10}]$ is at most $\frac{23}{15}<\frac{8}{5}<2$, the derivatives of order $k$ of $\tilde \o$ with respect to $\eta$ are multiplied by a factor smaller than $(\frac{5}{8})^k$. Consequently, we have
$$
M_1(\o)\le \frac{5}{8}(1+2M_1(\phi)\cdot \frac{\rho}{\kappa}).
$$
We can chose $\phi$ in such a way that $M_1(\phi)\le 2$. Then we get
$$
M_1(\o)\le \frac{5}{8}(1+4\frac{\rho}{\kappa})\le \frac{5}{8}(1+\frac{2}{5})=\frac{7}{8}<1.
$$
$$
M_k(\o)\le (\frac{5}{8})^k 2^k M_k(\phi)\cdot \frac{\rho}{\kappa^k}=(\frac{5}{4})^k M_k(\phi)\cdot \frac{\rho}{\kappa^k}:= \frac{D_k\rho}{\kappa^k}, \ \ k\ge 2.
$$
This completes the proof of Proposition \ref{prop:straight.curved}, with $D_k=(\frac{5}{4})^k M_k(\phi)$. $\square$

\section{Bounding $\nu_d(Z)$ via Metric entropy}\label{sec:entr.bounds.nu}
\setcounter{equation}{0}

In this section we show that $\bar \nu_d(Z)$, and hence $\nu_d(Z)$, can be bounded from above through metric entropy of $Z$. The idea is to mimic the approach of Integral Geometry on a discrete level. For a hypersurface $Y \subset {\mathbb R}^n$ one of the basic facts of Integral Geometry is that the average number of intersection points of $Y$ with a straight line in ${\mathbb R}^n$ is the $n-1$-area of $Y$. We will show that ``approximately the same'' is true for $Y$ replaced by any set $Z$ (in particular, finite), while the $n-1$-area of $Y$ is replaced with the appropriate covering numbers of $Z$.

\medskip

We use somewhat non-standard definition of covering numbers. We fix an orthonormal coordinate system in ${\mathbb R}^n$. For a given $\e>0$ we subdivide ${\mathbb R}^n$ into closed regular $\e$-cubes $Q_\alpha(\e)$, with the faces parallel to the coordinate hyperplanes, starting at the origin. The covering number $M(\e,Z)$ is then the number of the sub-cubes $Q_\alpha(\e)$, intersecting $Z$.

\medskip

Below we assume, as usual, that the set $Z$ is contained in the ball $\hat B^n$. Accordingly, in Theorem \ref{thm:int.geom.ent.1} below we consider only sub-cubes $Q_\alpha(\e)$ inside a $\frac{1}{10}$-neighborhood of $\hat B^n$. We will also always assume that $0\le \e \le \frac{1}{10\sqrt n}$. 

\medskip

Consider, as above, the set ${\cal L}(z_0)$ of all the straight lines $\ell$ in ${\mathbb R}^n$ passing through $z_0$ and intersecting the ball $\hat B^n$. We identify ${\cal L}(z_0)$ with the part of the unit sphere centered at $z_0$, and equipped with the Euclidean measure $\bar m_{n-1}$.

\medskip

Put $\theta_n=2^{3(n-1)}\frac{\bar \beta_{n-1}}{\beta_{n-1}},$ with $\bar \beta_{n-1}, \beta_{n-1}$ being the volume of ${\cal L}(z_0)$, and the volume of the $n-1$-dimensional unit ball, respectively.

\bt\label{thm:int.geom.ent.1}
Let $\e$, a point $z_0\in S^{n-1}$, and a natural $N$ be given. Then for each collection ${\cal U}$ of  $q\ge \theta_nN(\frac{1}{\e})^{n-1}$ sub-cubes $Q_\alpha(\e)$ inside $\hat B^n$, there exists a straight line $\ell$ in ${\mathbb R}^n,$ passing through $z_0$, which intersects at least $N$ of the sub-cubes $Q_\alpha(\e)$ in ${\cal U}$.
\et
\pr

Denote the $\e$-cubes $Q_\alpha(\e)$ in the collection ${\cal U}$ by $Q_j, \ j=1,\ldots,q,$ and define the functions $\psi_j(\ell)$ on ${\cal L}(z_0)$ as follows: $\psi_j(\ell)=1,$ if $\ell$ intersects the cube $Q_{j}$, and $\psi_j(\ell)=0$ otherwise. Then

\be\label{eq:ball.proj}
\int_{{\cal L}(z_0)} \psi_j(\ell) d\mu(\ell) \ge (\frac{1}{8})^{n-1}\e^{n-1}\beta_{n-1}.
\ee

Indeed, this integral is equal to the measure of the radial projection, from the point $z_0$, of the cube $Q_{j}$ to the unit sphere $S^{n-1}(z_0)$. The cube $Q_{j}$ contains the inscribed ball of radius $\frac{\e}{2}$. Hence, the radial projection, from the point $z_0$, of the cube $Q_{j}$ to the unit sphere $S^{n-1}(z_0)$, always contains a ball of radius at least $\frac{\e}{4}$, independently of the position of the cube inside $\hat B^n$, and of $z_0\in S^{n-1}$. Thus the bound (\ref{eq:ball.proj}) holds.

\medskip

Now we consider the function $\psi(\ell)=\sum_j \psi_j(\ell).$ Clearly, for each $\ell$, the function $\psi(\ell)$ is equal to the number of the cubes $Q_{j}$ that $\ell$ intersects. On the other hand, we have, by (\ref{eq:ball.proj}),
$$
\int_{{\cal L}(z_0)} \psi(\ell) d\mu(\ell)\ge q (\frac{1}{8})^{n-1}\e^{n-1}\beta_{n-1} \ge (\frac{1}{8})^{n-1}\beta_{n-1}\theta_n N  = \bar \beta_{n-1}N.
$$
This last integral is the average over $\ell \in {\cal L}(z_0)$ of the crossing number $\psi(\ell)$, and since the total measure of ${\cal L}(z_0)$ is equal to $\bar \beta_{n-1}$, we conclude that there exists some specific straight line $\bar \ell \in {\cal L}(z_0)$ for which $\psi(\bar\ell)\ge N.$ This completes the proof of Theorem \ref{thm:int.geom.ent.1}. $\square$.

\medskip
\medskip

In what follows we start with a given set $Z\subset \hat B^n$, and, for a given $\e$, consider the collection ${\cal U}$ of all the $\e$-sub-cubes intersecting with $Z$. We want to find a straight line $\ell$ in ${\mathbb R}^n$, passing through a given point $z_0$, and passing very close to certain $d+1$ points $z_1,\ldots,z_{d+1}$ in $Z$, such that their projections to $\ell$ are ``well-sepatated''. This will allow us, via Proposition \ref{prop:straight.curved}, to build a smooth curve $\o$ with ``small'' high-order derivatives, passing through $z_0$ and through $z_1,\ldots,z_{d+1}$.

\smallskip

Put $\xi_1=\xi_1(n,d)=2\theta_n n^2(d+1).$ As a consequence of Theorem \ref{thm:int.geom.ent.1} we obtain the following result:

\bp\label{prop:int.geom.ent.3}
Let a subset $Z\subset \hat B^n,$ a point $z_0\in S^{n-1}$, and positive numbers $\e$ and $\kappa>10\sqrt n\e$ be given. Assume that we have
$M(\e,Z)\ge \xi_1\kappa(\frac{1}{\e})^{n}.$ Then there exist $d+1$ points $z_1,\ldots,z_{d+1}$ in $Z$, and a straight line $\ell$ in ${\mathbb R}^n$, passing through $z_0$, such that the following conditions are satisfied:

\medskip

\noindent 1. The distance of each point $z_j$ to $\ell$ is at most $\sqrt n \e$.

\medskip

\noindent 2. The distance between the projections $\eta_j$ of the points $z_j$ to $\ell$ is at least $\kappa$.

\medskip

In particular, $\bar \nu_d(Z,z_0)\le \frac{\sqrt n \e}{\kappa^d}$.
\ep
\pr
Put $N=\frac{2n^2(d+1)\kappa}{\e}$, (which is significantly higher than $d+1$). Let ${\cal U}$ be the collection of all the $\e$-sub-cubes $Q_\alpha(\e)$, intersecting with $Z$, according to the definition of the covering number $M(Z,\e)$. Then, by the assumptions, the cardinality $|{\cal U}|$ satisfies the conditions of Theorem \ref{thm:int.geom.ent.1}. Indeed, we have
$$
|{\cal U}|=M(\e,Z)\ge \xi_1\kappa(\frac{1}{\e})^{n}=2\theta_n n^2(d+1)\kappa(\frac{1}{\e})^{n}= 
$$
$$
= \theta_n \frac{2n^2(d+1)\kappa}{\e}(\frac{1}{\e})^{n-1} = \theta_n N(\frac{1}{\e})^{n-1}.
$$
Theorem \ref{thm:int.geom.ent.1} now provides a straight line $\ell$ in ${\mathbb R}^n$, passing through $z_0$, and through at least $N$ \ $\e$-sub-cubes $Q_\alpha(\e)$ in ${\cal U}$. Take exactly $N$ of these cubes, and denote them $Q_j, \  j=1,\ldots,N.$ Next, in each of these cubes we fix a certain point $y_j\in Z\cap Q_j$. Each $y_j$ belongs to an $\e$-cube, which intersects $\ell$. Hence the distance of each point $y_j$ to $\ell$ is at most $\sqrt  n \e$.

\smallskip

Let $v=(v_1,\ldots,v_n)$ be the unit vector in the direction of $\ell$. We can assume that $|v_1|\ge \max_{i=2,\ldots,n} |v_i|.$ It is convenient to also assume that the first coordinates $\theta_j$ of $y_j$ are ordered: $\theta_j\le \theta_{j+1}, \ j=1,\ldots,N-1$.

\medskip

While some distances $\theta_{j+1}-\theta_{j}$ may be small (or zero), the following lemma shows that the average distance $\theta_{j+1}-\theta_{j}$ is of order at least $\frac{\e}{n}$:

\bl\label{lem:distance}
For each $j,s, \ j+s\le N,$ we have
$$
|\theta_{j+s}-\theta_j|\ge ([\frac{s}{n}]-1)\e.
$$
\el
\pr
Consider the ``$\e$-slices'' $\Sigma_p$ of ${\mathbb R}^n$:
$$
\Sigma_p=\{z=(x_1,\ldots,x_n)\in {\mathbb R}^n, \ p\e \le x_1 < (p+1)\e\}, \ p\in {\mathbb Z}.
$$
While crossing a certain slice $\Sigma_p$, and for each $i=2,\ldots,n$, the line $\ell$ may cross at most one of the coordinate hyperplanes $x_i=q\e$. Indeed, the first coordinate $v_1$ of $v$ is maximal in absolute value. Thus, as the coordinate $x_1$ of a point $w$ on $\ell$ goes from $p\e$ to $(p+1)\e$, the increment of each other coordinate of $w$ is at most one.

\smallskip

Therefore, $\ell$ can cross, inside the slide $\Sigma_p$, at most $n$ sub-cubes $Q_\alpha(\e)$. On the other hand, for given $j,s$, we have $y_j\in Q_j, \ y_{j+s}\in Q_{j+s}$, and between $Q_j$ and $Q_{j+s}$ the line $l$ crosses at least $s+1$ sub-cubes $Q_j,\ldots,Q_{j+s}$. By the comparison of the velocities of $\ell$ in the cooridnate directions, presented above, we conclude that while crossing these sub-cubes, the line $l$ must cross at least $[\frac{s}{n}]-1$ slices $\Sigma_p$. This completes the proof of Lemma \ref{lem:distance}. $\square$

\medskip

Now we put $s=[\frac{2n^2\kappa}{\e}]$ and pick the points $z_i$ to be $z_i=y_{si}, \ i=1,\ldots,d+1$. Notice that for $i=d+1$ we get
$si\le (d+1)\frac{2n^2\kappa}{\e}=N$, and thus we have enough points $y_j$ to pick all the required $z_i$.

\smallskip
\smallskip

Now, Lemma \ref{lem:distance} guarantees that the first coordinates $\eta_i$ of the points $z_i$ satisfy
$$
|\eta_{i+1}-\eta_i|\ge ([\frac{s}{n}]-1)\e \ge \frac{1.5n\kappa}{\e}\e= 1.5n\kappa,
$$
by the choice of $s$. Hence also the distances between the points $z_i$ satisfy $||z_{i+1}-z_{i}||\ge 1.5n\kappa.$ Finally, since the distance of each point $z_i$ to $\ell$ is at most $\sqrt  n \e$, we conclude, that the distances between the projections $\eta_i$ of the points $z_i$ to $\ell$ are at least $1.5n\kappa-2\sqrt  n \e > \kappa$. This completes the proof of Proposition \ref{prop:int.geom.ent.3}. $\square$

\medskip


\medskip

Now we have all the geometric tools we need in order to bound $\nu_d(Z)$ via the metric entropy of $Z$. Let $\xi=\xi_1(n,d)$ be the constant, defined in Proposition \ref{prop:int.geom.ent.3}. Put $\xi_2=\xi_2(n,d) = \sqrt  n \ \bar D(d)$, where $\bar D(d)$ is the constant, defined in Proposition \ref{prop:straight.curved}.

\bp\label{prop:int.geom.ent.4}
Let a subset $Z\subset \hat B^n$ be given. Assume that for certain $\e, \kappa, \ 0<\e < \frac{\kappa}{10\sqrt n},$ we have
$$
M(\e,Z)\ge \xi_1\kappa(\frac{1}{\e})^{n}.
$$
Then for each $z_0\in S^n$ there exists a curve $\o\in \O_d(Z,z_0)$ with
$$
\nu_d(\o)\le  \xi_2\cdot \frac{\e}{\kappa^{d+1}}.
$$
In particular, we have $\nu_d(Z)\le \xi_2\cdot \frac{\e}{\kappa^{d+1}}$.
\ep
\pr
We start with the straight line $\ell,$ passing through $z_0$, and with $d+1$ points $z_1,\ldots,z_{d+1}$ in $Z$, provided by Proposition \ref{prop:int.geom.ent.3}. Thus, the following conditions are satisfied:

\medskip

\noindent 1. The distance of each point $z_i$ to $\ell$ is at most $\sqrt  n \e$.

\medskip

\noindent 2. The minimal distance between the projections $\eta_i$ of the points $z_i$ to $\ell$ is at least $\kappa$.

\medskip

Now we apply Proposition \ref{prop:straight.curved}, and find a curve $\o\in \O_d(Z,z_0)$, passing through the points $z_1,\ldots,z_{d+1}$ and satisfying
$$
M_1(\o) \ \le \ 1, \ \ \ M_k(\o) \ \le \ \frac{D_k\sqrt  n \e}{\kappa^{k}}, \ \ \ k=2,3,\ldots ,
$$
We conclude that
$$
\nu_d(\o)=\max_{k=2}^{d+1} M_k(\o) \le \xi_2\cdot \frac{\e}{\kappa^{d+1}}, \  \  \xi_2 = \sqrt  n \ \max_{k=2}^{d+1} \ D_k=\sqrt  n \bar D(d).
$$
This completes the proof of Proposition \ref{prop:int.geom.ent.4}. $\square$

\section{Proof of main results}\label{Sec:proof.main.res}
\setcounter{equation}{0}

As a consequence of Proposition \ref{prop:int.geom.ent.4} we obtain the proof of our main result, which estimates the $d$-thickness $\nu_d(Z,z_0)$ in terms of the covering density of $Z$. It implies also the results of Theorems \ref{thm:main.intro1} and \ref{main:rigid.h.regular}, stated in the Introduction.

\medskip

For $\e>0$ we put $\kappa(\e)=\frac{M(\e,Z)\e^{n}}{\xi_1}$. This is the maximal possible $\kappa$, satisfying the condition                        $M(\e,Z)\ge\xi_1(n,d)\kappa(\frac{1}{\e})^{n}$ of Proposition \ref{prop:int.geom.ent.4}.

\medskip

Next we introduce the following geometric characteristic of sets $Z$:

\bd\label{def:d.density}
For $Z\subset \hat B^n$ we define $\zeta_d(Z)$ as
$$
\zeta_d(Z)=\sup_{\e: \ \frac{1}{10} \ \ge \ \e \ > \ 0, \ \kappa(\e) \ \ge \ 10\sqrt n \e} M(Z,e)\e^{n-\frac{1}{d+1}}.
$$
\ed
The $d$-density of $Z$ is a kind of a ``fractal volume'' of $Z$ in dimension $n-\frac{1}{d+1}$. Similar geometric characteristics appear in some other problems of fractal geometry (see e.g. \cite{Bru.Yom,Fri.Yom,Yom3}). We plan to present separately a more detailed their study.

\smallskip

Put $\xi_3:= \xi_2^{\frac{1}{d+1}}\xi_1$, where the constants $\xi_1, \ \xi_2$ were defined in Propositions \ref{prop:int.geom.ent.3} and \ref{prop:int.geom.ent.4}, respectively.

\bt\label{thm:main3}
For a subset $Z\subset \hat B^n$ we have
$$
\nu_d(Z)\le (\frac{\xi_3}{\zeta_d(Z)})^{d+1}.
$$
\et
\pr
Let $\e_0>0$ be the value of $\e$ for which the supremum in Definition \ref{def:d.density} is achieved. Put $\kappa_0=\kappa(\e_0)=\frac{M(\e_0,Z)\e_0^{n}}{\xi_1}$. By our definitions, the conditions of Proposition \ref{prop:int.geom.ent.4} are satisfied for $(\e,\kappa)=(\e_0,\kappa_0)$, and hence we have
$$
\nu_d(Z)\le \xi_2\cdot \frac{\e_0}{\kappa_0^{d+1}}=\frac{\xi_2\xi_1^{d+1}\e_0}{(M(\e_0,Z))^{d+1}\e_0^{n(d+1)}}.
$$
Extracting the $(d+1)$-st root from both parts of this inequality, we get
$$
(\nu_d(Z))^{\frac{1}{d+1}}\le \frac{\xi_2^{\frac{1}{d+1}}\xi_1\e_0^{\frac{1}{d+1}}}{M(\e_0,Z)\e_0^n} =\frac{\xi_3}{M(\e_0,Z)\e_0^{n-\frac{1}{d+1}}}=\frac{\xi_3}{\zeta_d(Z)}.
$$
Raising back to the $(d+1)$-st power, we finally get
$$
\nu_d(Z)\le (\frac{\xi_3}{\zeta_d(Z)})^{d+1}.
$$
This completes the proof of Theorem \ref{thm:main3}. $\square$

\medskip

Now we return to Theorems \ref{thm:main.intro1} and \ref{main:rigid.h.regular}, stated in the Introduction. Recall that the Minkowski (or ``box'', or ``entropy'') dimension $dim_e(Z)$ of $Z$ is defined as
$dim_e(Z)=\overline{\lim}_{\e \to 0}\frac{\log M(\e,Z)}{\log \frac{1}{\e}}$. As an immediate consequence of Theorem \ref{thm:main3} we obtain:

\bc\label{cor:ent.dim}
For $Z\subset \hat B^n,$ if $dim_e(Z)> n-\frac{1}{d+1},$ then $\nu_d(Z)=0$ and hence ${\cal RG}_d(Z)\ge C_4(n,d)>0$.
\ec
\pr
If $dim_e(Z)=n-\frac{1}{d+1}+\beta,$ with $\frac{1}{d+1} \ge \beta>0,$ then for $\e\to 0,$ we have
$$
M(\e,Z)\e^{n-\frac{1}{d+1}}\sim (\frac{1}{\e})^\beta\to \infty.
$$
In turn, $\kappa(\e)=\frac{M(\e,Z)\e^{n}}{\xi_1}\sim \e^{\frac{1}{d+1}-\beta}$. We conclude that for $\e\to 0$ we have $\frac{\kappa(\e)}{\e}\to \infty$. In particular, the condition $\kappa(\e) \ \ge \ 10\sqrt n \e$ of Definition \ref{def:d.density} is satisfied. Finally we get $\zeta(Z)=\infty$, and, by Theorem \ref{thm:main3}, $\nu_d(Z)=0$. Application of Theorem \ref{thm:rigidity.width} provides the required bound         ${\cal RG}_d(Z)\ge C_4(n,d)>0$ on the $d$-rigidity of $Z$. This completes the proof of Corollary \ref{cor:ent.dim} and of Theorem \ref{thm:main.intro1} (with $M=C_4(n,d)>0$).  $\square$

\medskip

Nest we consider finite sets $Z$ inside the cube $Q^n_s=[0,s]^n$ of size $s$. Let us recall that a set $Z\subset Q^n_s$ is called $h$-dense if in each $h$-sub-cube of $Q^n_s$ there is a point of $Z$. In particular, if we start with a regular $h/2$-grid in $Q^n_s$, and shift its points to distances at most $\frac{h}{8},$ we obtain an $h$-dense set.

\smallskip

Put $\xi=\frac{C_3}{2\xi_3^{d+1}}.$

\bt\label{main:rigid.h.dense}
Let $Z\subset Q^n_s$ be an $h$-dense set. Then if $h\le \xi s^{n(d+1)}$, then
$$
{\cal RG}_d(Z)\ge C_4(n,d)>0.
$$
\et
\pr
Put $\e=h$ in our definition of the covering number. By definition of an $h$-dense set we conclude that each $\e$-cube in the partition intersects $Z$. Therefore for $h\ll s$ we have $M(Z,h)= (1+o(h)) (\frac{s}{h})^n$. Hence
\be\label{eq:zeta}
\zeta_d(Z)\ge M(Z,h)h^{n-\frac{1}{d+1}}= (1+o(h))s^n h^{-\frac{1}{d+1}}.
\ee
Now assume that $h\le h_0=\xi s^{n(d+1)}$. Substituting into (\ref{eq:zeta}) we get
$$
\zeta_d(Z)\ge (1+o(h))s^n h^{-\frac{1}{d+1}}\ge (1+o(h_0))s^n h_0^{-\frac{1}{d+1}}=(1+o(h_0))\xi^{-\frac{1}{d+1}},
$$
and hence, by Theorem \ref{thm:main3},
$$
\nu_d(Z)\le (\frac{\xi_3}{\zeta_d(Z)})^{d+1}\le (1+o(h_0))\xi_3^{d+1}\xi < C_3(n,d),
$$
by our choice of $\xi$.

\medskip

Notice that $\kappa(h_0)=\frac{M(h_0,Z)h_0^{n}}{\xi_1}\sim \frac{s^n}{\xi_1}$. We conclude that 
$$
\kappa(h_0)\gg h_0\sim s^{n(d+1)}. 
$$
In particular, the condition $\kappa(h_0) \ \ge \ 10\sqrt n h_0$ of Definition \ref{def:d.density} is satisfied.

\medskip

Finally, by Theorem \ref{thm:rigidity.width}, we conclude that
$$
{\cal RG}_d(Z)\ge C_4(n,d)>0.
$$
This completes the proof of Theorems \ref{main:rigid.h.dense} and \ref{main:rigid.h.regular}. $\square$

\medskip

\bibliographystyle{amsplain}

\end{document}